\RequirePackage{fix-cm}
\let\origvec\vec
\documentclass[smallextended]{svjour3}       

\let\vec\origvec
\smartqed  
\usepackage [latin1]{inputenc}
\usepackage{amsmath}
\usepackage{mathrsfs}
\usepackage{amssymb}
\usepackage{graphicx}
\usepackage[justification=centering]{caption}
\usepackage{lineno}
\usepackage{color}
 \numberwithin{equation}{section}

\journalname{myjournal}
\begin{document}
\title{The global bifurcation of periodic internal waves with  point vortex and capillary effect}


\author{Guowei Dai \and Yong Zhang 
}


\institute{Guowei Dai \at
              School of Mathematical Sciences, Dalian University of Technology, Dalian, 116024, PR China \\
              \email{daiguowei@dlut.edu.cn.}           
              \and
            Yong Zhang (Corresponding author)\at
             School of Mathematical Sciences, Jiangsu University, Zhenjiang 212013, PR China \\
              \email{zhangyong@ujs.edu.cn}
}

\date{Received: date / Accepted: date}

\maketitle

\begin{abstract}
In this paper, we first construct two-dimensional periodic interface waves with point vortex and capillary effect and then obtain the global structure of the set of solutions.  This is done using the local and global bifurcation argument. Especially, we establish a global continuation theorem by using the degree for $C^1$ Fredholm mappings. As far as we know, the global result is new and general, which is promising to deal with other new phenomena in water waves.
\end{abstract}



\section{\bf Introduction}

\qquad In the past twenty years, rotational water waves have been completely considered. Indeed, the presence of wind forcing, temperature gradients, or even a slight heterogeneity in the density, would generate vorticity. The study of rotational steady waves essentially begins with Dubreil-Jacotinin 1934 (see \cite{Dubreil}), but the entire topic becomes hot until relatively recently. The main breakthrough was done by Constantin and Strauss in \cite{ConstantinS}, where they developed a systematic existence theory for two-dimensional periodic steady travelling gravity waves. Both the work of Constantin-Strauss and Dubreil-Jacotin begin with a simple observation: if there are no stagnation points in the flow, then one can use the stream function as a vertical coordinate to fix the domain. Doing so, one ultimately arrives at a quasilinear elliptic PDE on a strip (the quasilinearity coming from the change of coordinates) with a nonlinear boundary condition. Then one can build small-amplitude solutions by a perturbative argument based on shear flows. Then, using a degree theoretic continuation method, Constantin and Strauss were able to obtain finite-amplitude solutions. Following the publication of \cite{ConstantinS}, many authors have been able to generalize this approach, for instance \cite{Wahlen,Walsh1,Walsh2,Walsh3}.

However, we would like to emphasis that the degree theoretic continuation method used by Constantin and Strauss \cite{ConstantinS} was based on Rabinowitz's result, where the general Leray-Schauder degree theory developed by Kielhofer \cite{Kielhofer1} played a key role. In this paper, one of main aims is to extend the Rabinowitz's result to the Fredholm type operator equation. Since $C^1$ Fredholm mappings, in general, cannot be written as completely continuous vector fields, the Leray-Schauder degree theory and general Leray-Schauder degree theory cannot be used to show the desired conclusions.
Thus, we will use the degree for $C^1$ Fredholm mappings of index $1$, which is originally introduced in \cite{Pejsachowicz}. Then we obtain the following global continuation theorem.
\begin{theorem}
Let $X$ and $Y$ be real Banach spaces and $F:\mathbb{R}\times X\longrightarrow Y$ be a $C^1$ map.
Suppose that

1. $D_{u}F(\lambda, u)$ is a Fredholm operator with index $0$ for all $(\lambda,u)\in \mathbb{R}\times X$;

2. $F$ is proper on any closed bounded subsets of $\mathbb{R}\times X$;

3. $u=u_0$ is the unique solution of $F\left(\lambda_0,u\right)=0$;

4. $D_u F\left(\lambda_0,u_0\right)\in L(X,Y)$ is bijective.


\noindent Then $\mathscr{S}$ contains a pair of unbounded components $\mathcal{C}^+$, $\mathcal{C}^-$
in $\mathbb{R}_{\lambda_0}^+\times X$, $\mathbb{R}_{\lambda_0}^-\times X$ respectively and $\mathcal{C}^+\cap\mathcal{C}^-=\left\{\left(\lambda_0,u_0\right)\right\}$, where $\mathbb{R}_{\lambda_0}^+=\left[\lambda_0,+\infty\right)$ and $\mathbb{R}_{\lambda_0}^-=\left(-\infty,\lambda_0\right]$.
\end{theorem}

As an application of Theorem 1, then we consider the classical problem of establishing the existence of a free boundary dividing two superposed incompressible, inviscid, and immiscible fluids under the influence
of gravity and capillary tension. This situation arises in countless applications, with a particularly important example being internal waves propagating along a pycnocline or thermocline in the ocean.
Recent years have seen enormous progress made in understanding the dynamics of the interface in internal water waves. These waves are usually created by the presence of two different layers of water combined with a certain configuration of current. A well-known example is the wave-current interactions in the equatorial Pacific Ocean \cite{ConstantinB,ConstantinI}, where a thermocline (the sharp boundary between warm and deeper cold water) would occur which propagates beneath the surface, along a density interface with a characteristic speed of $0.1m/s$.  As far as we know, the first small-amplitude interfacial periodic traveling capillary-gravity waves on finite depth was constructed by Amick and Turner \cite{AmickT}, where the vorticity of fluid was ignored. Recently, the work \cite{WalshOS} allowed for the global vorticity in the interfacial wind waves by using the semi-hodograph transformation. We also want to mention \cite{AmbroseSW}, where the new coordinates were chosen and the Rabinowitz's global bifurcation theorem was applied to prove the existence of pure capillary interfacial waves and gravity-capillary interfacial waves, where the fluid interface may be not a graph over the horizontal. In addition, we refer to \cite{ChuDE,Sinambela,Tony,Wheeler} for more information on two-layer or many-layer density stratified water waves of global vorticity.

Unlike in previous studies,
here we consider solutions which are periodic perturbations of a background shear flow with localised vorticity rather than global vorticity (see the following Fig. 1).
\begin{figure}[ht]
\centering
\includegraphics[width=0.75\textwidth]{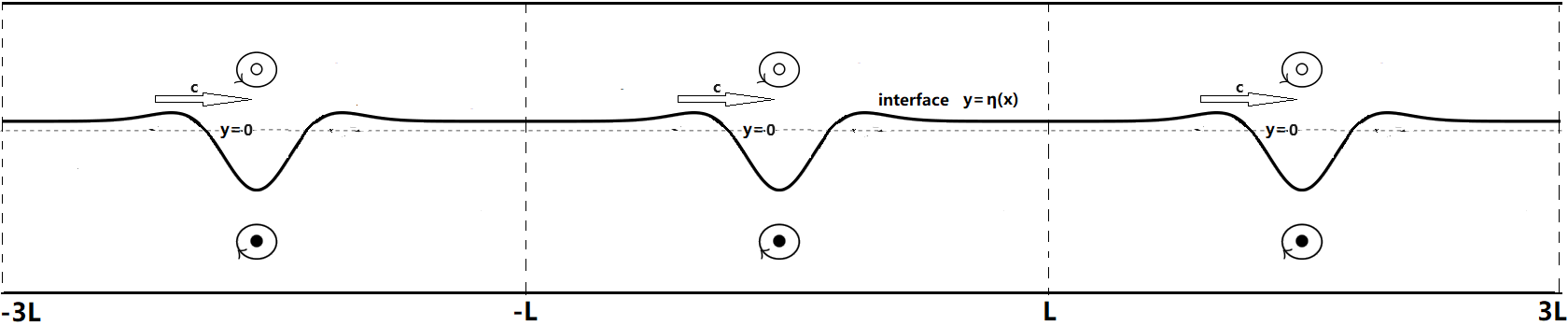}
\caption{The periodic internal waves with point vortices}
\label{fig1}
\end{figure}

Now let us state the other two main results in this paper.
\begin{theorem} \label{thm1.1} (Local bifurcation)
Let $s\geq 3$, then there is an open interval $I$ including $0$ and a $C^{\infty}$ curve $\mathcal{C}_{loc}$
$$
I\quad \rightarrow \quad \mathcal{O}^s\times \mathbb{R}
$$
$$
\varepsilon \quad \rightarrow \quad (\eta(\varepsilon), \bar{\xi}(\varepsilon) ,\xi(\varepsilon) ,c(\varepsilon) ,\varepsilon)
$$
of solutions to the Zakharov-Criag-Sulem system (\ref{eq2.16})-(\ref{eq2.18}), where $|\varepsilon|<\varepsilon_0$, $\mathcal{O}^s\subset H^s_{even}(L,-L)\times H^s_{even}(L,-L)\times H^s_{even}(L,-L)\times \mathbb{R}$ and
$$
(\eta(0), \bar{\xi}(0) ,\xi(0) ,c(0) ,0)=(0, 0, 0, 0).
$$
Moreover, $\mathcal{C}_{loc}$ comprises all solutions of (\ref{eq2.16})-(\ref{eq2.18}) in a sufficiently small neighbourhood of the origin point in $\mathcal{O}^s\times \mathbb{R}$.
\end{theorem}

\begin{theorem} \label{thm1.2} (Global bifurcation)
There exits a connected set $\mathcal{C}\subset \mathcal{O}^s\times \mathbb{R}$ of solutions to the Zakharov-Criag-Sulem system (\ref{eq2.16})-(\ref{eq2.18}) with $\mathcal{C}_{loc}\subset\mathcal{C}$. Moreover, one of the following alternatives must hold:

(i) there is a sequence $\{ (\eta_n,\bar{\xi}_n,\xi_n,c_n,  \varepsilon_n) \}\subset \mathcal{C}$ which is unbounded in $\mathcal{O}^s\times \mathbb{R}$, or


(ii) along some sequence in $\{ (\eta_n,\bar{\xi}_n,\xi_n,c_n, \varepsilon_n) \}\subset \mathcal{C}$ with $|\eta_n(0)|=d$.
\end{theorem}

The rest of this paper is arranged as follows. We begin Section 2 by introducing the governing equations for two-dimensional steady periodic interfacial waves with point vortex and rewrite the government equation into Zakharov-Craig-Sulem formulation. In Section 3, we establish the existence result of small-amplitude interfacial waves with point vortex by classic implicit function theorem. The Section 4 is devoted to establishing a new global implicit function theorem. In Section 5, we apply the global implicit function theorem to extend the local curve of solutions obtained in Section 3 to the global.

\section{\bf Statement of the problem}
It's known that the vertical stratification of the ocean is a common phenomenon, where a shallow layer of relatively warm water (superlayer of less dense) overlies a much deeper layer of cold water (underlayer of higher density). The interface between the warm water and the cold water is called the thermocline. Without loss of generality, here we only consider steady travelling internal waves of the thermocline or interface with point vortex in a period $[-L,L]$ for some $L>0$. (see Fig. 2).
\begin{figure}[ht]
\centering
\includegraphics[width=0.75\textwidth]{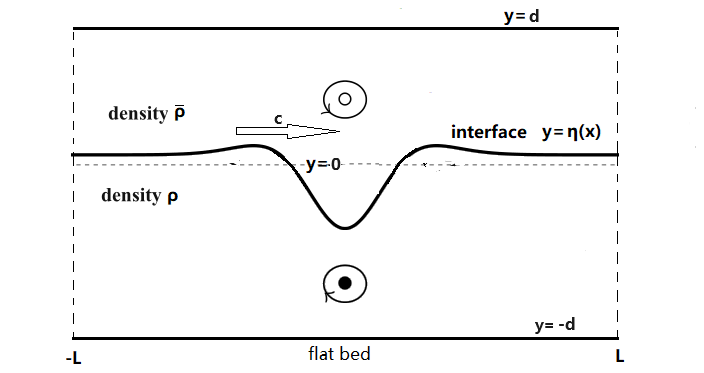}
\caption{The internal wave in one period}
\label{fig2}
\end{figure}

\subsection{\bf The governing equations}
Now we shall present the governing equations for two-dimensional internal waves traveling at the boundary region between two fluids with different densities under the rigid lid condition at the top. The bottom of the underlayer is assumed to be flat, denoted by $y=-d$. Let $y=d$ be the upper boundary of the superlayer and $y=\eta(t, x)$ be the interface oscillating about the horizontal line $y=0$. Then it is obvious that the mean depth of superlayer and underlayer is $d$ with $d>0$ being some constant. We also suppose that the internal waves move with a constant wavespeed $c>0$. Now let us denote the superlayer by $\Omega_{\eta}^a:=\{(x,y):-L<x<L,~~ \eta(t,x)<y<d\}$ and denote the underlayer by $\Omega_{\eta}^b:=\{(x,y):-L<x<L,~~ -d<y<\eta(t,x)\}$, then the governing equations can be expressed by the following nonlinear problem (Euler equations)
\begin{equation}\label{eq2.1}
\left\{\begin{array}{llll}
{\bar{u}_{x}+\bar{v}_{y}=0} &{~\text{in}~\Omega_{\eta}^a},\\
{\bar{u}_t+\bar{u}\bar{u}_{x}+\bar{v}\bar{u}_{y}=-\frac{1}{\bar{\rho}}\bar{P}_{x}} &{~\text{in}~\Omega_{\eta}^a},\\
{\bar{v}_t+\bar{u}\bar{v}_{x}+\bar{v}\bar{v}_{y}=-g-\frac{1}{\bar{\rho}}\bar{P}_{y}} &{~\text{in}~\Omega_{\eta}^a},
\end{array}\right.
\end{equation}
and
\begin{equation}\label{eq2.2}
\left\{\begin{array}{llll}
{u_{x}+v_{y}=0} &{~\text{in}~\Omega_{\eta}^b},\\
{u_t+uu_{x}+vu_{y}=-\frac{1}{\rho}P_{x}} &{~\text{in}~\Omega_{\eta}^b},\\
{v_t+uv_{x}+vv_{y}=-g-\frac{1}{\rho}P_{y}} &{~\text{in}~\Omega_{\eta}^b},
\end{array}\right.
\end{equation}
subjected to the boundary conditions
\begin{equation}\label{eq2.3}
\left\{\begin{array}{llll}
{\bar{v}=0} &{~\text{on}~y=d},\\
{v=0} &{~\text{on}~y=-d},
\end{array}\right.
\left\{\begin{array}{llll}
{P=\bar{P}-\sigma\frac{\eta_{xx}}{(1+\eta_x^2)^\frac{3}{2}}} &{~\text{on}~y=\eta(t,x)},\\
{v=\eta_t+u\eta_x} &{~\text{on}~y=\eta(t,x)},\\
{\bar{v}=\eta_t+\bar{u}\eta_x} &{~\text{on}~y=\eta(t,x)},
\end{array}\right.
\end{equation}
where $(\bar{u},\bar{v})$ are the velocity field, $\bar{P}$ is the pressure, $\bar{\rho}$ is the water's density of the superlayer $\Omega_{\eta}^a$ and we denote by $(u,v)$, $P$ and $\rho$ the velocity, pressure and density of the colder water in the underlayer $\Omega_{\eta}^b$ beneath the interface with $\sigma>0$ being the surface tension coefficient and $g$ the gravitational acceleration.

To be more precise, we are interested in finding solutions of problem (\ref{eq2.1})-(\ref{eq2.3}) with "almost irrotation", that is the vorticity consists of point vortex at some place. By a point vortex we mean that the vorticity is given by a $\delta$-function at this place. In fact, the waves with point vortex are the simplest form of waves with compactly supported vorticity.

\subsection{\bf The point vorticity and travelling waves}
By taking the curl of Euler equations in (\ref{eq2.1}) and (\ref{eq2.2}), we can obtain that
\begin{equation}\label{eq2.4}
\bar{\omega}_t+((\bar{u},\bar{v})\cdot\nabla)\bar{\omega}=0,\quad  \omega_t+((u,v)\cdot\nabla)\omega=0
\end{equation}
which imply that the vorticity $\bar{\omega}$ and $\omega$ are transported by the vector filed $(\bar{u},\bar{v})$ and $(u,v)$, respectively. This fact admits us to expect that if the vorticity is a point vortex at some time, then it will be always a point vortex in future and be transported by the flow. However, it should be emphasized that, for now, the point vortex is too singular to give even a weak solution to the Euler equations. Completely faithful adherence to (\ref{eq2.4}) would require that the center $\mathbf{z}=\mathbf{z}(t)\subset\Omega_{\eta}^b$ (or $\bar{\mathbf{z}}=\bar{\mathbf{z}}(t)\subset\Omega_{\eta}^a$) of the point vortex in underlayer (or superlayer) move controlled by $\frac{d\mathbf{z}}{dt}=\mathbf{U}(t,\mathbf{z})=(u,v)$ (or $\frac{d\bar{\mathbf{z}}}{dt}=\overline{\mathbf{U}}(t,\mathbf{z})=(\bar{u},\bar{v})$). But if $\omega$ (or $\bar{\omega}$) is a Dirac measure, the magnitude of the velocity diverges as one approaches $\mathbf{z}$ (or $\bar{\mathbf{z}}$), which makes the ODE ill-posed. Helmholtz and Kirchhoff independently arrived at a resolution to this issue. The key physical reasoning is that the point vortex is not transported by the full velocity field because it does not self-advect; it is driven only by the irrotational part of $\mathbf{U}$ at $\mathbf{z}$ (or $\overline{\mathbf{U}}$ at $\bar{\mathbf{z}}$). To make sense of this, observe that we can always split the velocity field as
\begin{equation}
\mathbf{U}(x,y,t)=\mathbf{V}(x,y,t)+\mathbf{W}(x,y,t), \quad  \overline{\mathbf{U}}(x,y,t)=\overline{\mathbf{V}}(x,y,t)+\overline{\mathbf{W}}(x,y,t)  \nonumber
\end{equation}
where
$$\nabla \times \mathbf{V}=0,\quad \nabla \times \mathbf{W}=\omega, \quad \nabla \cdot \mathbf{V}=0,\quad \nabla \cdot \mathbf{W}=0$$
and
$$\nabla \times \overline{\mathbf{V}}=0, \quad \nabla \times \overline{\mathbf{W}}=\bar{\omega}, \quad \nabla \cdot \overline{\mathbf{V}}=0, \quad \nabla \cdot \overline{\mathbf{W}}=0.$$

Based on this splitting, we may take $\mathbf{W}$ (or $\overline{\mathbf{W}}$) as the part of the velocity field $\mathbf{U}$ (or $\overline{\mathbf{U}}$) generated by the point vortex, is singular and divergence free around the point $\mathbf{z}(t)=(x_0(t), y_0(t))$ in $\Omega_{\eta}^b$ (or $\bar{\mathbf{z}}(t)=(\bar{x}_0(t), \bar{y}_0(t))$ in $\Omega_{\eta}^a$). Based on the assumed curl and divergence of these vector fields above, it's not difficult to find that there are stream functions $\phi$ (or $\bar{\phi}$) and $\varphi$ (or $\bar{\varphi}$), determined uniquely modulo constants by $\mathbf{V}$ (or $\overline{\mathbf{V}}$) and $\mathbf{W}$ (or $\overline{\mathbf{W}}$) such that
$$
\mathbf{V}=\nabla^{\bot}\phi, \quad \mathbf{W}=\nabla^{\bot}\varphi \quad \text{and} \quad \overline{\mathbf{V}}=\nabla^{\bot}\bar{\phi}, \quad \overline{\mathbf{W}}=\nabla^{\bot}\bar{\varphi},
$$
where $\nabla^{\bot}=(-\partial_y,\partial_x)$ with $\phi$ and $\bar{\phi}$ being harmonic and
\begin{equation}\label{eq2.5}
\Delta\varphi=\delta_{\mathbf{z}}, \quad \Delta\bar{\varphi}=\delta_{\mathbf{z}}.
\end{equation}

In the following, we will suppose that $\omega=\varepsilon\delta_{\mathbf{z}}$  and $\bar{\omega}=-\varepsilon\delta_{\bar{\mathbf{z}}}$ with $\varepsilon=\int_{\Omega^b_{\eta}}\omega dX$ being strength of point vortex $\omega$. It is worth noting that the point vortex $\bar{\omega}$ in superlayer is caused by the point vortex $\omega$ in underlayer. Thus we can think of $\bar{\omega}$ as a phantom vortex that has the same strength (see Fig. 1).
Then it's easy to see that
\begin{equation}\label{eq2.6}
\mathbf{U}=\nabla^{\bot}\phi+\varepsilon\nabla^{\bot}\varphi,\quad \overline{\mathbf{U}}=\nabla^{\bot}\bar{\phi}+\varepsilon\nabla^{\bot}\bar{\varphi}.
\end{equation}
This decomposition is unique once appropriate boundary conditions have been prescribed.
It follows from (\ref{eq2.5}) that
\begin{equation}\label{eq2.7}
\varphi=\Gamma(\cdot-\mathbf{z})+\varphi_{\mathcal{H}},\quad \bar{\varphi}=\Gamma(\cdot-\bar{\mathbf{z}})+\bar{\varphi}_{\mathcal{H}},
\end{equation}
where $\Gamma:=\frac{1}{4\pi}\log(x^2+y^2)$ is the fundamental solution of $\Delta$ in $\mathbb{R}^2$ and $\varphi_{\mathcal{H}}$ and $\bar{\varphi}_{\mathcal{H}}$ are two harmonic functions. Based on arguments on point vortex above, we will take $\varphi_{\mathcal{H}}=-\Gamma(\cdot-\bar{\mathbf{z}})$ and $\bar{\varphi}_{\mathcal{H}}=-\Gamma(\cdot-\mathbf{z})$ in the following. Since the $\bar{\omega}$ is thought of a phantom vortex, here we only consider the vorticity function in underlayer.
The Helmholtz-Kirchhoff model states that the point vortex moves should satisfy
\begin{equation}\label{eq2.8}
\frac{\partial \mathbf{z}}{\partial t}=(\nabla^{\bot} \phi+\varepsilon\nabla^{\bot} \varphi_{\mathcal{H}})(\mathbf{z}).
\end{equation}
Indeed, if the singular self-interaction term $\varepsilon\nabla^{\bot}\Gamma(\cdot-\mathbf{z})$ were removed, the vector fields are analytic (see \cite{ChuDE}).

For the travelling waves, it is convenient to introduce a new reference frame moving along the wave by $(x-ct,y)\mapsto (x,y)$. Then the flow is time independent and it is natural to assume that $\eta(x,t)=\eta(x-ct)$ and
$$
\overline{\mathbf{U}}(x,y,t)=\mathbf{\overline{U}}(x-ct,y), \quad \mathbf{U}(x,y,t)=\mathbf{U}(x-ct,y),
$$
$$
\overline{P}(x,y,t)=\overline{P}(x-ct,y), \quad P(x,y,t)=P(x-ct,y)
$$
for all relevant $x$, $y$ and $t$, where positive and negative $c$ is corresponding to waves moving in the positive and negative $x$-directions, respectively.
In the new steady variables, the governing equations (\ref{eq2.1})-(\ref{eq2.3}) will become
\begin{equation}\label{eq2.9}
\left\{\begin{array}{llll}
{\bar{u}_{x}+\bar{v}_{y}=0} &{~\text{in}~\Omega_{\eta}^a},\\
{(\bar{u}-c)\bar{u}_{x}+\bar{v}\bar{u}_{y}=-\frac{1}{\bar{\rho}}\bar{P}_{x}} &{~\text{in}~\Omega_{\eta}^a},\\
{(\bar{u}-c)\bar{v}_{x}+\bar{v}\bar{v}_{y}=-g-\frac{1}{\bar{\rho}}\bar{P}_{y}} &{~\text{in}~\Omega_{\eta}^a},
\end{array}\right.
\end{equation}
and
\begin{equation}\label{eq2.10}
\left\{\begin{array}{llll}
{u_{x}+v_{y}=0} &{~\text{in}~\Omega_{\eta}^b},\\
{(u-c)u_{x}+vu_{y}=-\frac{1}{\rho}P_{x}} &{~\text{in}~\Omega_{\eta}^b},\\
{(u-c)v_{x}+vv_{y}=-g-\frac{1}{\rho}P_{y}} &{~\text{in}~\Omega_{\eta}^b},
\end{array}\right.
\end{equation}
subjected to the boundary conditions
\begin{equation}\label{eq2.11}
\left\{\begin{array}{llll}
{\bar{v}=0} &{~\text{on}~y=d},\\
{v=0} &{~\text{on}~y=-d},
\end{array}\right.
\left\{\begin{array}{llll}
{P=\bar{P}-\sigma\frac{\eta_{xx}}{(1+\eta_x^2)^\frac{3}{2}}} &{~\text{on}~y=\eta(t,x)},\\
{v=(u-c)\eta_x} &{~\text{on}~y=\eta(t,x)},\\
{\bar{v}=(\bar{u}-c)\eta_x} &{~\text{on}~y=\eta(t,x)}
\end{array}\right.
\end{equation}
and the point vorticity equation given in (\ref{eq2.8}) would reduce to
\begin{equation}\label{eq2.12}
(-\phi_y-\varepsilon(\varphi_{\mathcal{H}})_y, \phi_x+\varepsilon(\varphi_{\mathcal{H}})_x)(\mathbf{z})=(c,0).
\end{equation}

\subsection{\bf The Zakharov-Craig-Sulem formulation}

In this subsection, we will reduce the problem (\ref{eq2.9})-(\ref{eq2.12}) to an one-dimensional one on the interface by using the Zakharov-Craig-Sulem formulation, which was first introduced by Zakharov in \cite{Zakharov}, and then later put on a strict mathematical basis by Criag and Sulem in \cite{CraigS} and \cite{CraigSS}. The original Zakharov-Craig-Sulem formulation depends heavily on the fluid being irrotational, but this holds strictly near the thermocline or interface $y=\eta(x)$ in this paper. This is why we require that the point vortices are far from the interface $y=\eta(x)$ in this paper.

We can use the momentum of conversation in (\ref{eq2.9}) and (\ref{eq2.10}) to deduce that
\begin{equation}\label{eq2.13}
\left\{\begin{array}{llll}
{\bar{\rho}\nabla\left( -c\bar{u}+\frac{1}{2}|\overline{\mathbf{U}}|^2+gy \right)+\nabla \bar{P}=0,} &{~\text{on}~y=\eta(x)},\\
{\rho\nabla\left( -cu+\frac{1}{2}|\mathbf{U}|^2+gy \right)+\nabla P=0,} &{~\text{on}~y=\eta(x)},
\end{array}\right.
\end{equation}
where we use the fact of point vortices far from the interface $y=\eta(x)$. It follows from (\ref{eq2.13}) that
\begin{equation}
\begin{array}{llll}
{c(\rho u-\bar{\rho}\bar{u})+\frac{1}{2}\bar{\rho}|\overline{\mathbf{U}}|^2-\frac{1}{2}\rho|\mathbf{U}|^2+(\bar{\rho}-\rho)gy+\bar{P}-P=Q,} &{~\text{on}~y=\eta(x)}, \nonumber
\end{array}
\end{equation}
where $Q$ is a real constant.
Considering the dynamic boundary condition in (\ref{eq2.11}) on $y=\eta(x)$, we have that
\begin{equation}\label{eq2.14}
\begin{array}{llll}
{c(\rho u-\bar{\rho}\bar{u})+\frac{1}{2}\bar{\rho}|\overline{\mathbf{U}}|^2-\frac{1}{2}\rho|\mathbf{U}|^2+(\bar{\rho}-\rho)gy+\sigma\frac{\eta_{xx}}{\left(1+\eta_x^2\right)^{\frac{3}{2}}}=Q}
\end{array}
\end{equation}
on $y=\eta(x)$.
On the other hand, it follows from the conversation of mass in (\ref{eq2.9}) and (\ref{eq2.10}) and the kinematic rigid boundary conditions in (\ref{eq2.11}) that
\begin{equation}\label{eq2.15}
\left\{\begin{array}{llll}
{\Delta\bar{\phi}=0} &{~\text{in}~\Omega_{\eta}^a},\\
{\bar{\phi}_{x}=0} &{~\text{on}~y=d},\\
{\bar{\phi}=\bar{\xi}(x)} &{~\text{on}~y=\eta(x)},
\end{array}\right.\quad
\left\{\begin{array}{llll}
{\Delta\phi=0} &{~\text{in}~\Omega_{\eta}^b},\\
{\phi_{x}=0} &{~\text{on}~y=-d},\\
{\phi=\xi(x)} &{~\text{on}~y=\eta(x)},
\end{array}\right.
\end{equation}
where $\bar{\xi}(x)$ and $\xi(x)$ are called the trace of stream function $\bar{\phi}(x,y)$ and $\phi(x,y)$ respectively limited on the interface $y=\eta(x)$. We would like to mention that if knowing $\bar{\phi}$ and $\phi$ , one will recover the velocity field $\overline{\mathbf{U}}$ or $\mathbf{U}$ and the pressure $\bar{P}$ or $P$
through (\ref{eq2.9})-(\ref{eq2.11}). Considering the elliptic problems (\ref{eq2.15}), we are thus led to find a set of equations that determines $\eta$, $\bar{\xi}$ and $\xi$. To this end, it is convenient to introduce the Harmonic extension operator
$$
\bar{H}(\eta): \bar{\xi}\mapsto \bar{\phi}, \quad H(\eta): \xi\mapsto \phi
$$
and the Dirichlet-Neumann operator
$$
\bar{G}(\eta): \bar{\xi}\mapsto \sqrt{1+\eta_x^2}\frac{\partial \bar{\phi}}{\partial \nu}\mid_{y=\eta(x)}, \quad G(\eta): \xi\mapsto \sqrt{1+\eta_x^2}\frac{\partial \phi}{\partial \nu}\mid_{y=\eta(x)},
$$
where $\bar{\phi}$ and $\phi$ solve (\ref{eq2.15}) and $\nu$ is the unit outward normal on $y=\eta(x)$. Now we will use these operators to reformulate the equation (\ref{eq2.14}), the kinematic interfacial boundary conditions in (\ref{eq2.11}) and vorticity equation (\ref{eq2.12}) in a way that only involves $(\eta, \bar{\xi}, \xi, c)$.

Note that
$$
\bar{\xi}_{x}(x)=\bar{\phi}_{x}+\eta_x\bar{\phi}_{y},\quad \xi_x=\phi_{x}+\eta_x\phi_{y},
$$
$$
\bar{G}(\eta)\bar{\xi}=\eta_x\bar{\phi}_{x}-\bar{\phi}_{y},\quad G(\eta)\xi=\phi_{y}-\eta_x\phi_{x},
$$
where the right-hand side is evaluated on $y=\eta(x)$. Inverting these relations and inserting them in (\ref{eq2.14}) and the kinematic interfacial boundary conditions in (\ref{eq2.11}), we can obtain that
\begin{eqnarray} \label{eq2.16}
~& &c\left(\bar{\rho}\frac{\bar{\xi}_{x}\eta_x+\bar{G}(\eta)\bar{\xi}}{1+\eta_x^2}+\varepsilon\bar{\rho}\bar{\varphi}_{y}-\rho\frac{\xi_{x}\eta_x+G(\eta)\xi}{1+\eta_x^2}-\varepsilon\rho\varphi_{y}\right)
\nonumber\\
& &+ \frac{\bar{\rho}}{2}\left(\left(\frac{\bar{G}(\eta)\bar{\xi}+\eta_x\bar{\xi}_x}{1+\eta_x^2}+\varepsilon\bar{\varphi}_y \right)^2+\left( \frac{\bar{\xi}_x-\eta_x\bar{G}(\eta)\bar{\xi}}{1+\eta_x^2}+\varepsilon\bar{\varphi}_x\right)^2\right) \nonumber\\
& &-\frac{\rho}{2}\left(\left(\frac{G(\eta)\xi+\eta_x\xi_x}{1+\eta_x^2}+\varepsilon\varphi_y \right)^2+\left( \frac{\xi_x-\eta_xG(\eta)\xi}{1+\eta_x^2}+\varepsilon\varphi_x\right)^2\right) \nonumber\\
& &+(\bar{\rho}-\rho)g\eta+\sigma\frac{\eta_{xx}}{\left(1+\eta_x^2\right)^{\frac{3}{2}}}-Q=0
\end{eqnarray}
and
\begin{equation}\label{eq2.17}
\left\{\begin{array}{llll}
{\xi_x+\varepsilon\varphi'+c\eta_x=0},\\
{\bar{\xi}_x+\varepsilon\bar{\varphi}'+c\eta_x=0},
\end{array}\right.
\end{equation}
where $\varphi'=\frac{d\varphi(x,\eta(x))}{dx}$ and $\bar{\varphi}'=\frac{d\bar{\varphi}(x,\eta(x))}{dx}$.
On the other hand, we have some point vortices carried by a interfacial solitary capillary-gravity wave with these point vortices being far from the boundary. We impose the following ansatz
$$
\varphi_{\mathcal{H}}=- \Gamma(\cdot-\mathbf{\bar{z}})\quad \text{and} \quad \bar{\varphi}_{\mathcal{H}}=- \Gamma(\cdot-\mathbf{z}),
$$
which are harmonic inside $\Omega_{\eta}^a$ and $\Omega_{\eta}^b$ respectively.
Therefore, it follows from discussion above and (\ref{eq2.12}) that
\begin{eqnarray} \label{eq2.18}
c=-(H(\eta)\xi)_y(\mathbf{z})+c_1\varepsilon,
\end{eqnarray}
where $c_1=\Gamma_{y}(\mathbf{z}-\mathbf{\bar{z}})$.
Now, we will focus on the governing equations in Zakharov-Craig-Sulem formulation (\ref{eq2.16})-(\ref{eq2.18}). It is not difficult find that the unknown functions would be $(\eta, \bar{\xi}, \xi, c)$ in Zakharov-Craig-Sulem formulation and the strength of point vortex $\varepsilon$ will be chosen as a parameter.

\section{Small-amplitude solitary waves with point vortex by implicit function theorem}
To prove the existence of a family of solitary interfacial waves of small amplitude with point vortices, now let us introduce a suitable functional analytic frame. For convenience, we define the spaces by
$$
X^s:=H^s_{even}(-L,L)\times H^s_{even}(-L,L)\times H^s_{even}(-L,L)\times \mathbb{R},
$$
$$
Y^s:=H^{s-2}_{even}(-L,L)\times H^{s}_{even}(-L,L)\times H^{s}_{even}(-L,L)\times \mathbb{R},
$$
where $H^s_{even}(-L,L):=\{ f\in H^s(-L,L): f~ \text{is even}\}$, $s\geq3$ and set
$$
\mathcal{O}^s:=\{(\eta,\bar{\xi},\xi,c)\in X^s : d(\bar{\mathbf{z}},\mathcal{S})>0, ~ d(\mathbf{z},\mathcal{S})>0\},
$$
which means the center of point vortex $\bar{\mathbf{z}}$ and $\mathbf{z}$ are far from the interface $\mathcal{S}:=\{(x,y) |y=\eta(x)\}$.
We proceed to introduce five maps that together will form the basis for our argument. For $s\geq3$, we define $F_1: \mathcal{O}^s\times \mathbb{R}\rightarrow H^{s-2}_{even}(-L,L)$ by
\begin{eqnarray}
F_1(\eta,\bar{\xi}, \xi, c, \varepsilon)& =&c\left(\bar{\rho}\frac{\bar{\xi}_{x}\eta_x+\bar{G}(\eta)\bar{\xi}}{1+\eta_x^2}+\varepsilon\bar{\rho}\bar{\varphi}_{y}-\rho\frac{\xi_{x}\eta_x+G(\eta)\xi}{1+\eta_x^2}-\varepsilon\rho\varphi_{y}\right)
\nonumber\\
& &+ \frac{\bar{\rho}}{2}\left(\left(\frac{\bar{G}(\eta)\bar{\xi}+\eta_x\bar{\xi}_x}{1+\eta_x^2}+\varepsilon\bar{\varphi}_y \right)^2+\left( \frac{\bar{\xi}_x-\eta_x\bar{G}(\eta)\bar{\xi}}{1+\eta_x^2}+\varepsilon\bar{\varphi}_x\right)^2\right) \nonumber\\
& &-\frac{\rho}{2}\left(\left(\frac{G(\eta)\xi+\eta_x\xi_x}{1+\eta_x^2}+\varepsilon\varphi_y \right)^2+\left( \frac{\xi_x-\eta_xG(\eta)\xi}{1+\eta_x^2}+\varepsilon\varphi_x\right)^2\right) \nonumber\\
& &+(\bar{\rho}-\rho)g\eta+\sigma\frac{\eta_{xx}}{\left(1+\eta_x^2\right)^{\frac{3}{2}}}-Q, \nonumber
\end{eqnarray}
the maps $F_2: \mathcal{O}^s\times \mathbb{R}\rightarrow H^{s}_{even}(-L,L)$ and $F_3: \mathcal{O}^s\times \mathbb{R}\rightarrow H^{s}_{even}(-L,L)$ by
\begin{eqnarray}
F_2((\eta,\bar{\xi}, \xi, c, \varepsilon)=\bar{\xi}+\varepsilon\bar{\varphi}+c\eta, \nonumber
\end{eqnarray}
\begin{eqnarray}
F_3(\eta,\bar{\xi}, \xi, c, \varepsilon)=\xi+\varepsilon\varphi+c\eta\nonumber
\end{eqnarray}
and the maps $F_4: \mathcal{O}^s\times \mathbb{R}\rightarrow \mathbb{R}$ by
\begin{eqnarray}
F_4(\eta,\bar{\xi},\xi,c,\bar{\varepsilon}, \varepsilon)=c+(H(\eta)\xi)_y(\mathbf{z})-c_1\varepsilon. \nonumber
\end{eqnarray}

We can now define $F:\mathcal{O}^s\times \mathbb{R}\rightarrow Y^s$ by
$$
F:=(F_1,F_2,F_3,F_4),
$$
then we will find solutions of the equation
\begin{eqnarray}\label{eq3.1}
F(\eta,\bar{\xi}, \xi,c, \varepsilon)=0.
\end{eqnarray}
It is easy to check that $F(0,0,0,0,0)=0$, which means that the origin is a trivial solution. Next we will turn out that in a small neighborhood of the origin in $\mathcal{O}^s\times \mathbb{R}$ there is a unique curve of nontrivial solutions parametrized by the vortex strength parameter $\varepsilon$.\\
~\\
{\bf Proof of Theorem 2:} As mentioned above, the origin is a trivial solution. In the following, we will apply the implicit function theorem to finish the proof. To this end, let us compute the first partial derivatives of $F$ at the origin $(0,0,0,0,0)$. A direct calculation yields that
\begin{equation}  \label{eq3.2}
D_XF(0,0,0,0,0)=\left[
\begin{array}{cccc}
(\bar{\rho}-\rho)g+\sigma\partial_{xx}  \qquad&  0 \qquad& 0 \qquad& 0\\
~ & ~\\
0  \qquad&  I \qquad& 0 \qquad& 0\\
~ & ~\\
0  \qquad&  0 \qquad& I \qquad& 0\\
~ & ~\\
0    \qquad& 0 \qquad& \left(H(0)\cdot\right)_y(\mathbf{z}) \qquad& 1\\
\end{array}
\right], \nonumber
\end{equation}
where the subscript $X$ denotes the partial derivative with respect to the variable $(\eta, \bar{\xi}, \xi,c)$ in $\mathcal{O}^s$.

Now, we claim the operator $D_XF(0,0,0,0,0)\in L(\mathcal{O}^s\times\mathbb{R}, Y^s)$ is an isomorphism. Indeed, the operator
$$
(\bar{\rho}-\rho)g+\sigma\partial_{xx} : H^s_{even}(-L,L)\rightarrow  H^{s-2}_{even}(-L,L)
$$
corresponds to the Fourier multiplier $(\bar{\rho}-\rho)g-\sigma\zeta^2$. Since $(\bar{\rho}-\rho)g<0$ and $\sigma>0$, this operator is invertible.
The other three operators on the diagonal are identity operators, and therefore trivially invertible. Thus, we conclude that the $D_XF(0,0,0,0,0)\in L(\mathcal{O}^s\times\mathbb{R}, Y^s)$ is an isomorphism.
Then we can use the implicit function theorem to conclude the result. \qed

\section{ \bf The proof of Theorem 1}
In this section, we mainly establish a new global implicit function theorem, i.e Theorem 1. Note that Kielhofer \cite{Kielhofer} also obtained a similar result. However, our conditions and conclusions here are different from \cite[Theorem II.6.1]{Kielhofer}. Here we require $F$ is weaker than the one where $F$ is required to be $C^2$ and our conclusion has only one choice, which is better than his alternative conclusion. It is worth noting that we require $u=u_0$ being the unique solution of $F\left(\lambda_0,u\right)=0$ which is stronger than his.

 Let $X$ and $Y$ be real Banach spaces. We shall investigate the structure of the set of nontrivial solutions for the following nonlinear operator equation
\begin{equation}\label{operater}
F(\lambda,u)=0,\,\, (\lambda,u)\in \mathbb{R}\times X,
\end{equation}
where $F:\mathbb{R}\times X\longrightarrow Y$
is $C^1$ smooth. We assume that $u=u_0$ is the unique solution of $F\left(\lambda_0,u\right)=0$ for some $\lambda_0\in \mathbb{R}$ and $D_u F\left(\lambda_0,u_0\right):X\longrightarrow Y$ is Fredholm map of index $0$.
Clearly, $\left(\lambda_0,u_0\right)$ is a solution of $F(\lambda,u)=0$, which is called trivial solution.
Let $\mathscr{S}$ denote the closure of the set of nontrivial solutions of (\ref{operater}) which contains $\left(\lambda_0,u_0\right)$.

When $Y=X$ and $F(\lambda,u)=u-G(\lambda,u)$ with $G:\mathbb{R}\times X\rightarrow X$ is compact and continuous.
Rabinowitz \cite[Theorem 3.2]{Rabinowitz} established a global continuation theorem for (\ref{operater}) as follows.

\begin{proposition}
If $G\left(\lambda_0,u\right)\equiv u_0$, then $\mathscr{S}$ contains a pair of unbounded components $\mathcal{C}^+$, $\mathcal{C}^-$
in $\mathbb{R}_{\lambda_0}^+\times X$, $\mathbb{R}_{\lambda_0}^-\times X$ respectively and $\mathcal{C}^+\cap\mathcal{C}^-=\left\{\left(\lambda_0,u_0\right)\right\}$, where $\mathbb{R}_{\lambda_0}^+=\left[\lambda_0,+\infty\right)$ and $\mathbb{R}_{\lambda_0}^-=\left(-\infty,\lambda_0\right]$.
\end{proposition}

Now let us extend the Rabinowitz's result to the Fredholm type operator equation.
Let $\Phi_0(X,Y)$, $GL(X,Y)$ and $K(X,Y)$ be the set of linear Fredholm operators of index $0$, linear invertible operators and linear compact operators, respectively.
Given a continuous path $\alpha:[a,b]\longrightarrow \Phi_0(X,Y)$, it has been proved in \cite{Fitzpatrick0} that there exists a continuous path $\beta:[a,b]\longrightarrow GL(X,Y)$ such that
\begin{equation}
\beta(\lambda)\alpha(\lambda)-I\in K(X,X).\nonumber
\end{equation}
The path $\beta$ is called a \emph{parametrix} of $\alpha$.
Set $I-\beta(\lambda)\alpha(\lambda):=P(\lambda)$.
If $\alpha(a)$ and $\alpha(b)$ are isomorphisms, the \emph{parity} (see \cite{FitzpatrickP}\cite{FitzpatrickPejsachowicz}\cite{FitzpatrickRabier}) of $\alpha$ on $[a,b]$ is defined by
\begin{equation}
\sigma(\alpha,[a,b])=\deg(\beta(a)\alpha(a))\deg(\beta(b)\alpha(b)),\nonumber
\end{equation}
where $\deg(\beta(a)\alpha(a))$ is the Leray-Schauder degree of $(I-P(a))(x)=0$ with respect to any open set containing $0$.

A map $G:X\longrightarrow Y$ is called \emph{proper} if the inverse
image of each compact subset of $Y$ is a compact subset of $X$. Let $\Omega$ be an open subset. A Fredholm mapping $G\in C^1(X,Y)$ of index $1$ (equivalently, $D_uG(\cdot)$ is
Fredholm operator with index $0$) such that
$G\big|_{\overline{\Omega}}$ is proper will be called $\Omega$-\emph{admissible}. We call that $p$ is a \emph{base point} of $G$
if $D_u G(p)\in GL(X,Y)$. If $G$ is $\Omega$-admissible and $y\in Y\setminus G(\partial\Omega)$ is a regular value of $G\big|_{{\Omega}}$, then $G^{-1}\{y\}\cap\Omega$ is finite:
\begin{equation}
G^{-1}\{y\}\cap\Omega:=\left\{x_1,\ldots,x_k\right\}.\nonumber
\end{equation}
If $G^{-1}\{y\}\cap\Omega\neq\emptyset$ and $p$ is a base point of $G$ corresponding to $y$.
Let $\gamma_i\in C^0\left(\left[a_i,b_i\right],X\right)$ be any curve joining $p$ and $x_i$, where $i\in\{1,\ldots,k\}$.
Define the \emph{degree}
\begin{equation}
d_p(G,\Omega,y)=\Sigma_{i=1}^k\sigma_i,\nonumber
\end{equation}
where $\sigma_i=\sigma\left(\left(D_uG\right)\circ\gamma_i,\left[a_i,b_i\right]\right)$.
If $G^{-1}\{y\}\cap\Omega=\emptyset$, we set $d_p(G,\Omega,y):=0$. If $y$ is singular value, define $d_p(G,\Omega,y)$ by $d_p(G,\Omega,y):=d_p(G,\Omega,z)$ where $z\in Y\setminus G(\partial\Omega)$ is
any regular value of $G\big|_\Omega$ lying in a sufficiently small neighborhood of $y$ in $Y$ (the existence of such regular values is ensured by the Quinn-Sard theorem).
From the definition, we clearly see that
\begin{equation}\label{translation invariance}
d_p(G-y,\Omega,0)=d_p(G,\Omega,y).
\end{equation}
In particular, we find that
\begin{equation}\label{normativeness}
d_p(I-y,\Omega,0)=1.
\end{equation}
Indeed, taking $\gamma=I$, $b=a$ and $\beta=I$, we see that
\begin{equation}
\sigma(\left(D_uG\right)\circ\gamma,\left[a,b\right])=\deg(\beta(a)\alpha(a))\deg(\beta(b)\alpha(b))=\deg(I)\deg(I)=1,\nonumber
\end{equation}
which implies (\ref{normativeness}).\\

To show Theorem 1, we recall the following famous topological lemma due to Whyburn \cite[Theorem 9.3 of Chapter I]{Whyburn}.
\begin{lemma}
Let $K$ be a compact metric space and $A$ and $B$ disjoint closed subsets of $K$. Then either there exists a sub-continuum of $K$ meeting both $A$
and $B$ or $K=K_A\cup K_B$, where $K_A$, $K_B$ are disjoint compact subsets of $K$ containing $A$ and $B$, respectively.
\end{lemma}

Since $u=u_0$ is the unique solution of $F\left(\lambda_0,u\right)=0$, the local constancy of \cite[Corally 5.3]{Pejsachowicz} and normalization of \cite[Corally 5.5]{Pejsachowicz} implies that (\ref{operater}) possesses a solution for all $\lambda$ such that $\left\vert \lambda-\lambda_0\right\vert$ sufficiently small.
Let $\mathcal{C}^+$, $\mathcal{C}^-$ denote respectively the components of $\mathscr{S}$ in $\mathbb{R}_{\lambda_0}^+\times X$, $\mathbb{R}_{\lambda_0}^-\times X$.\\

Based on the above topological lemma, we can establish the following crucial result.
\begin{lemma}
Under the assumptions of Theorem 3, suppose
that there does not exist an unbounded sub-continuum of $\mathscr{S}\cup\left\{\left(\lambda_0,u_0\right)\right\}$ which meets $\left(\lambda_0,u_0\right)$.
Then there exists a bounded open subset $U$ of $\mathcal{O}$ containing $\mathcal{C}^+\cup \mathcal{C}^-$ such that
$\partial U\cap \mathscr{S}=\emptyset$.
\end{lemma}
\textbf{Proof.} Letting $\mathcal{C}=\mathcal{C}^+\cup \mathcal{C}^-$, then $\mathcal{C}$ is bounded. Since $F$ is proper o any closed bounded subsets of $\mathbb{R}\times X$,
$\mathcal{C}$ is also compact. Let $\mathscr{N}$ be a bounded open neighborhood of $\mathcal{C}$.
Set $K:=\overline{\mathscr{N}}\cap \mathscr{S}$. Since $\mathscr{S}$ is locally compact in $\mathbb{X}:=\mathbb{R}\times X$,
$K$ is a compact metric space under the induced topology from $\mathbb{X}$.
Let $A=\mathcal{C}$ and $B=\partial \mathscr{N}\cap \mathscr{S}$. By the construction of $\mathscr{N}$, we see that
$A\cap \partial \mathscr{N}=\emptyset$.
Using Lemma 1, there exist two disjoint compact subsets $K_A$ and $K_B$ containing $A$ and $B$, respectively, such that $K=K_A\cup K_B$.
Let $U$ be any $\epsilon$ neighborhood of $K_A$ in $\mathbb{X}$ where $\epsilon$ is less than the distance between $K_A$ and $K_B$.
Then $\partial U\cap \mathscr{S}=\emptyset$.\qed
~\\

Now we can present the argument of Theorem 1.\\
\\
\textbf{Proof of Theorem 1.}  If there does not exist unbounded branch $\mathcal{C}^\nu$, by Lemma 2, there exists a bounded open subset $U^\nu=\left\{(\lambda,u)\in U:\nu\left(\lambda-\lambda_0\right)\geq0\right\}$ containing $\mathcal{C}^\nu$ such that
$\partial U\cap \mathscr{S}=\emptyset$.
Since $D_u F\left(\lambda_0,u_0\right):X\longrightarrow Y$ is Fredholm map of index $0$ and $U^\nu$ is connected,
$D_u F\left(\lambda,u\right): X\longrightarrow Y$ is Fredholm map of index $0$ for any $(\lambda,u)\in U^\nu$.

Since $F$ is $C^1$ and $D_u F\left(\lambda_0,u_0\right)\in L(X,Y)$ is bijective, $p=0$ is a base point for both $F(\lambda,\cdot)$ for $\lambda$ near $\lambda_0$.
So $d_0\left(F(\lambda,\cdot),U_\lambda^\nu,0\right)$ is well defined for $\lambda$ near $\lambda_0$.
In particular,
since $\Phi\left(\lambda_0,u\right)=u-u_0$, using (\ref{translation invariance}), we have that
\begin{equation}
c=d_0\left(u,U_{\lambda_0}^\nu,u_0\right).\nonumber
\end{equation}
Further, by (\ref{normativeness}), we obtain that $c=1$.

We use $P\left(\overline{U}\right)$ to denote the projection of $\overline{U}$ on the axis of parameter $\lambda$.
Let $\lambda_*=\inf \left\{\lambda:\lambda\in P\left(\overline{U}\right)\right\}$ and $\lambda^*=\sup \left\{\lambda:\lambda\in P\left(\overline{U}\right)\right\}$.
We see that $U_{\lambda_*}=\emptyset=U_{\lambda^*}$ and
\begin{equation}
d_0\left(F(\lambda,\cdot),U_{\lambda_*},0\right)=d_0\left(\Phi(\lambda,\cdot),U_{\lambda^*},0\right)=0.\nonumber
\end{equation}
Note that $F^{-1}(\lambda,\cdot)\{0\}\cap U_\tau=F^{-1}(\lambda,\cdot)\{0\}\cap U_\lambda$ for $\vert \lambda-\tau\vert$ small enough for all $\lambda\in  P\left(\overline{U}\right)$. If $p_\lambda$ is a base point for $F(\lambda,\cdot)$, then $p_\lambda$ remains a base point for $F(\tau,\cdot)$ with $\vert \tau-\lambda\vert$ small enough.
Next, we will only consider the case of $\nu=+$ because the argument for $\nu=-$ is similar.
From the compactness of $\left[\lambda_0,\lambda^*\right]$, we find an integer $m>0$ such that, setting $\lambda_i=\lambda_0+i\left(\lambda^*-\lambda_0\right)/m$, $i\in\{0,\ldots,m\}$, a common base point $p_i\in X$
exists for $F(\lambda,\cdot)$, $\lambda\in\left[\lambda_i,\lambda_{i+1}\right]$, $i\in\{0,\ldots,m-1\}$, $p_0=p_m=0$.

According to the excision property of degree \cite[Theorem 5.2]{Pejsachowicz}, we obtain
\begin{equation}
d_{p_i}\left(F\left(\lambda,\cdot\right),U_{\lambda},0\right)=d_{p_i}\left(F\left(\lambda,\cdot\right),U_{\lambda_i},0\right),\,\,
\lambda\in\left[\lambda_i,\lambda_{i+1}\right], i\in\{0,\ldots,m-1\}.\nonumber
\end{equation}
It follows that
\begin{equation}
d_{p_i}\left(F\left(\lambda_{i+1},\cdot\right),U_{\lambda_{i+1}},0\right)=d_{p_i}\left(F\left(\lambda_{i+1},\cdot\right),U_{\lambda_i},0\right),\,\,
i\in\{0,\ldots,m-1\}.\nonumber
\end{equation}
By the homotopy invariance of the degree \cite[Theorem 4.1]{Pejsachowicz}, we obtain
\begin{equation}
d_{p_i}\left(F\left(\lambda_{i},\cdot\right),U_{\lambda_{i}},0\right)=\varepsilon d_{p_i}\left(F\left(\lambda_{i+1},\cdot\right),U_{\lambda_i},0\right),\,\,
i\in\{0,\ldots,m-1\},\nonumber
\end{equation}
where $\varepsilon$ is the parity of $D_uF(t\lambda_i+(1-t)\lambda_{i+1},\cdot)$, $t\in[0,1]$.

By invariance under reparametrizations \cite[Theorem 3.5]{Pejsachowicz} and \cite[Theorem 3.1]{Pejsachowicz}, we find that
\begin{equation}
\sigma\left(D_u F\left(t\lambda_i+(1-t)\lambda_{i+1},p_i\right),[0,1]\right)=\sigma\left(D_u F(\lambda,p_i),\left[\lambda_i,\lambda_{i+1}\right]\right)=1.\nonumber
\end{equation}
Thus, we conclude that
\begin{equation}
d_{p_i}\left(F\left(\lambda_{i},\cdot\right),U_{\lambda_{i}},0\right)=d_{p_i}\left(F\left(\lambda_{i+1},\cdot\right),U_{\lambda_{i+1}},0\right),\,\,
i\in\{0,\ldots,m-1\}.\nonumber
\end{equation}
Passing from the base point $p_i$ to the base point $p_{i+1}$ in the right hand side of the above relation, we obtain that
\begin{equation}
d_{p_i}\left(F\left(\lambda_{i},\cdot\right),U_{\lambda_{i}},0\right)=\varepsilon_i d_{p_{i+1}}\left(F\left(\lambda_{i+1},\cdot\right),U_{\lambda_{i+1}},0\right),\,\,
i\in\{0,\ldots,m-1\},\nonumber
\end{equation}
where $\varepsilon_i$ is the parity of $D_uF\left(\lambda_{i+1},\gamma_i\right)$ and $\gamma_i$ is a continuous curve in $X$ joining $p_i$ to $p_{i+1}$.

Let $\Gamma_i$ be the curve in $\mathbb{R}\times X$ joining $\left(\lambda_i,p_i\right)$ to $\left(\lambda_{i+1},p_{i+1}\right)$
and consisting of the line segment $\left(\lambda,p_i\right)$, $\lambda\in\left[\lambda_i,\lambda_{i+1}\right]$ followed by $\left(\lambda_{i+1},\gamma_i\right)$.
By \cite[Theorem 3.4]{Pejsachowicz}, we find that
\begin{eqnarray}
\sigma\left(D_uF\circ\Gamma_i\right)& =&\sigma\left(D_u F(\lambda,p_i),\left[\lambda_i,\lambda_{i+1}\right]\right)\sigma\left(D_uF\left(\lambda_{i+1},\gamma_i\right)\right)
\nonumber\\
&= &\varepsilon_i\sigma\left(D_u F(\lambda,p_i),\left[\lambda_i,\lambda_{i+1}\right]\right). \nonumber
\end{eqnarray}
Since $\lambda_m=\lambda^*$ and $p_0=p_m=0$, the curve $\Gamma$ in $\mathbb{R}\times X$
obtained from the following curves $\Gamma_0$, $\Gamma_1$, $\ldots$, $\Gamma_{m-1}$, successively, is a continuous curve joining
$\left(\lambda_0,0\right)$ to $\left(\lambda^*,0\right)$.
By \cite[Theorem 3.4]{Pejsachowicz}, the parity $\varepsilon$ of $D_uF\circ\Gamma$ equals the product $\varepsilon_0\cdots \varepsilon_{m-1}$.
Therefore, we obtain that
\begin{equation}
1=d_{0}\left(F\left(\lambda_0,\cdot\right),U_{\lambda_0},0\right)=\varepsilon_0\cdots \varepsilon_{m-1} d_{0}\left(F\left(\lambda^*,\cdot\right),U_{\lambda^*},0\right)=0.\nonumber
\end{equation}
This contradiction implies the desired conclusions.\qed

\section{The global bifurcation}
In this section, we will use the new global implicit function Theorem 1 to obtain the global bifurcation result Theorem 3. Based on the content of Theorem 1 and the arguments in Section 3, it is enough for us to check the operator $F$ is locally proper and the linearized operator $D_{X}F(X,\varepsilon)$ is Fredholm with 0-index for $(X,\varepsilon)=(\eta,\bar{\xi},\xi,c,\varepsilon)\in \mathcal{O}^s\times \mathbb{R}$.

\begin{lemma} (Proper)
Let $\mathcal{U}\subset Y$ be compact and $\mathbf{K}\subset \overline{\mathcal{O}^s}\times \mathbb{R}$ be closed and bounded. Then
$\mathbf{K}\cap F^{-1}(\mathcal{U})$ is compact in $X\times \mathbb{R}$.
\end{lemma}
\textbf{Proof.}
Let $\mathcal{U}$ and $\mathbf{K}$ be given as above. Define $\mathcal{V}:=\mathbf{K}\cap F^{-1}(\mathcal{U})\subset \mathcal{O}^s\times \mathbb{R}$ and let $v_n:=\{ (\eta_n,\bar{\xi}_n,\xi_n,c_n,\varepsilon_n) \}$ be a bounded sequence in $\mathcal{V}$. By construction
$$
u_n=(f_n,\bar{g}_n,g_n,h_n) :=F(\eta_n,\bar{\xi}_n,\xi_n,c_n,\varepsilon_n)\subset \mathcal{U}
$$
is bounded in $Y$. It follows from the compactness of $\mathcal{U}$ that there is a sequence $\{u_n\}$ converges to some $u$ in $Y$. Now we need to construct a convergent subsequence of $\{v_n\}$ in $\mathcal{V}$.

Since $c_n,\varepsilon_n \in \mathbb{R}$ is bounded, then it is easy to see that $c_n\rightarrow c$ and $\varepsilon_n\rightarrow\varepsilon$ along a subsequence for some $c,\varepsilon\in \mathbb{R}$. In addition, it follows from the definitions of operator equation (\ref{eq3.1}) and $u_n,v_n$ that
\begin{eqnarray}
f_n& =&c_{n}\left(\bar{\rho}\frac{\bar{\xi}'_{n}\eta'_n+\bar{G}(\eta_{n})\bar{\xi}_{n}}{1+\eta'^2_n}+\varepsilon_{n}\bar{\rho}\bar{\varphi}_{y}-
\rho\frac{\xi'_{n}\eta'_n+G(\eta_{n})\xi_{n}}{1+\eta'^2_n}-\varepsilon_{n}\rho\varphi_{y}\right)
\nonumber\\
& &+ \frac{\bar{\rho}}{2}\left(\left(\frac{\bar{G}(\eta_{n})\bar{\xi}_{n}+\eta'_n\bar{\xi}'_n}{1+\eta'^2_n}+\varepsilon_{n}\bar{\varphi}_y \right)^2+\left( \frac{\bar{\xi}'_n-\eta'_n\bar{G}(\eta_{n})\bar{\xi}_{n}}{1+\eta'^2_n}+\varepsilon_{n}\bar{\varphi}_x\right)^2\right) \nonumber\\
& &-\frac{\rho}{2}\left(\left(\frac{G(\eta_{n})\xi_{n}+\eta'_n\xi'_n}{1+\eta'^2_n}+\varepsilon_{n}\varphi_y \right)^2+\left( \frac{\xi'_n-\eta'_nG(\eta_{n})\xi_{n}}{1+\eta'^2_n}+\varepsilon_{n}\varphi_x\right)^2\right) \nonumber\\
& &+(\bar{\rho}-\rho)g\eta_{n}+\sigma\frac{\eta''_{n}}{\left(1+\eta'^2_n\right)^{\frac{3}{2}}}-Q. \nonumber
\end{eqnarray}
For convenience, let us rewrite
$$
f_n:= M_n(\eta_n,\bar{\xi}_n,\xi_n,c_n,\varepsilon_n)+\sigma\frac{\eta''_{n}}{\left(1+\eta'^2_n\right)^{\frac{3}{2}}}
$$
and denote
$$
M_{nm}:=M_n-M_m,\qquad \kappa_{nm}=\sigma\frac{\eta''_{n}}{\left(1+\eta'^2_n\right)^{\frac{3}{2}}}-\sigma\frac{\eta''_{m}}{\left(1+\eta'^2_m\right)^{\frac{3}{2}}}.
$$
Then it is easy to see that
\begin{eqnarray}\label{eq4.4}
f_{nm}:=f_n-f_m=M_{nm}+\kappa_{nm}.
\end{eqnarray}
It follows from the convergency of $\{f_n\}$ in $H^{s-2}(-L,L)$ by assumption that
\begin{eqnarray}\label{eq4.5}
\partial_{x}^{s-2}f_{nm}\rightarrow 0, \quad in \quad L^2(-L,L).
\end{eqnarray}
Since $\{M_n\}$ is bounded in $H^{s-1}(-L,L)$, the sequence $\{\partial_{x}^{s-2}M_{nm}\}$ is bounded in $H^1(-L,L)$, which means that
$\{\partial_{x}^{s-2}M_{nm}\}$ is bounded in $C^{0,\frac{1}{2}}(-L,L)$ by imbedding theorem. Hence the sequence $\{\partial_{x}^{s-2}M_{nm}\}$ is precompact in $C^{0,\alpha}(-L,L)$ for any $\alpha\in [0,\frac{1}{2})$. Then we can extract a convergent subsequence in $C^{0,\alpha}(-L,L)$ such that it is convergent in $L^2(-L,L)$, that is to say
\begin{eqnarray}\label{eq4.6}
\partial_{x}^{s-2}M_{nm}\rightarrow 0, \quad in \quad L^2(-L,L).
\end{eqnarray}

It follows from (\ref{eq4.4})-(\ref{eq4.6}) that
$$
\partial_{x}^{s-2}\kappa_{nm}\rightarrow 0, \quad in \quad L^2(-L,L).
$$
Since $\{\eta'_n\}$ is bounded in $H^{s-1}(-L,L)$, we can extract a subsequence such that $\{\left(1+\eta'^2_n\right)^\frac{1}{2}\}$ converges in $H^k(-L,L)$ for $k<s-1$. For $s\geq 3$, thus we have that
\begin{eqnarray}
\|\eta''_n-\eta''_m\|_{H^{s-2}}& =&\sigma\|\kappa_n\left(1+\eta'^2_n\right)^\frac{3}{2}-\kappa_m\left(1+\eta'^2_m\right)^\frac{3}{2}\|_{H^{s-2}}
\nonumber\\
& &\leq \sigma\|\kappa_{nm}\|_{H^{s-2}}\|\left(1+\eta'^2_n\right)^\frac{3}{2}\|_{H^{s-2}} \nonumber \\
& &+ \sigma\|\kappa_m\|_{H^{s-2}}\|\left(1+\eta'^2_n\right)^\frac{3}{2}-\left(1+\eta'^2_m\right)^\frac{3}{2}\|_{H^{s-2}}\rightarrow 0, \nonumber
\end{eqnarray}
which implies that $\{\eta_n\}$ has a convergent subsequence in $H^s(-L,L)$.

Similarly, we have that
$$
\bar{g}_n=\bar{\xi}_n+\varepsilon_n\bar{\varphi}+c_n\eta_n, \quad g_n=\xi_n+\varepsilon_n\varphi+c_n\eta_n.
$$
Therefore, it is obvious that $\{\bar{\xi}_n\}$ and $\{\xi_n\}$ have a convergent subsequence in $H^s(-L,L)$ due to the compactness of $\{\bar{g}_n\}$, $\{g_n\}$, $\{\varepsilon_n\}$, $\{c_n\}$ and $\{\eta_n\}$.
\qed
~\\

From now on, we will denote by $D_XF$ the following Fr\'{e}chet gradient
$$
D_XF:=(D_{\eta}F,~ D_{\bar{\xi}}F, ~D_{\xi}F, ~D_{c}F)
$$
and we would proceed to prove the remaining Fredholm property.

\begin{lemma} (Fredholm)
At each $(\eta,\bar{\xi},\xi,c,\varepsilon)\in \mathcal{O}^s\times \mathbb{R}$, the linearized operator $D_XF(\eta,\bar{\xi},\xi,c,\varepsilon)$ is a Fredholm operator of index 0.
\end{lemma}
\textbf{Proof.}
Fix $(\eta,\bar{\xi},\xi, c)\in \mathcal{O}^s$, we will establish the following estimate: for all $(\eta^*,\bar{\xi}^*,\xi^*, c^*)\in X^s$
\begin{equation}\label{eq4.7}
\|(\eta^*,\bar{\xi}^*,\xi^*, c^*)\|_{X^s}\leq \|D_XF(\eta,\bar{\xi},\xi,c,\varepsilon)(\eta^*,\bar{\xi}^*,\xi^*, c^*)\|_{Y^s}+\|(\eta^*,\bar{\xi}^*,\xi^*, c^*)\|_{Z^s},
\end{equation}
where $Z^s$ is a space such that $X^s$ is compactly embedded into $Z^s$ with constant may depend on $(\eta,\bar{\xi},\xi,c,\varepsilon)$.
Once the estimate (\ref{eq4.7}) is proved, we have that $F$ is semi-Fredholm. By the connectedness of $\mathcal{O}^s$ and the fact that the point of bifurcation $D_XF$ is an isomorphism, then we can deduce that $D_XF(\eta,\bar{\xi},\xi,c,\varepsilon)$ is a Fredholm operator of index 0.

Thus, let us establish the estimate (\ref{eq4.7}) by using the elliptic regularity. By a simple computation, we have that
\begin{eqnarray}\label{eq4.8}
~ & &
F_{1\eta}\eta^*\nonumber\\
& &=c\bar{\rho}\frac{(\bar{\xi}_x\eta_x^*+\bar{G}_{\eta}(\eta^*)\bar{\xi})(1+\eta_x^2)-2(\bar{\xi}_x\eta_x+\bar{G}(\eta)\bar{\xi})\eta_x\eta^*_x}{\left(1+\eta_x^2\right)^2}
\nonumber\\
& &-c\rho\frac{(\xi_x\eta_x^*+G_{\eta}(\eta^*)\xi)(1+\eta_x^2)-2(\xi_x\eta_x+G(\eta)\xi)\eta_x\eta^*_x}{\left(1+\eta_x^2\right)^2} \nonumber\\
& &+ \bar{\rho}\left(\frac{\bar{G}(\eta)\bar{\xi}+\eta_x\bar{\xi}_x}{1+\eta_x^2}+\varepsilon\bar{\varphi}_y \right)\frac{\left(\bar{G}_{\eta}(\eta^*)\bar{\xi}+\eta_x^*\bar{\xi}_x\right)\left(1+\eta_x^2\right)-2\left(\bar{G}(\eta)\bar{\xi}+\eta_x\bar{\xi}_x\right)\eta_x\eta_x^*}{\left(1+\eta_x^2\right)^2} \nonumber\\
& &+ \bar{\rho}\left( \frac{\bar{\xi}_x-\eta_x\bar{G}(\eta)\bar{\xi}}{1+\eta_x^2}+\varepsilon\bar{\varphi}_x\right)
\frac{\left(-\eta_x^*\bar{G}(\eta)\xi-\eta_x\bar{G}_{\eta}(\eta^*)\bar{\xi}\right)\left(1+\eta_x^2\right)-2\left(\bar{\xi}_x-\eta_x\bar{G}(\eta)\bar{\xi}\right)\eta_x\eta_x^*}{\left(1+\eta_x^2\right)^2} \nonumber\\
& &-\rho\left(\frac{G(\eta)\xi+\eta_x\xi_x}{1+\eta_x^2}+\varepsilon\varphi_y \right)\frac{\left(G_{\eta}(\eta^*)\xi+\eta_x^*\xi_x\right)\left(1+\eta_x^2\right)-2\left(G(\eta)\xi+\eta_x\xi_x\right)\eta_x\eta_x^*}{\left(1+\eta_x^2\right)^2} \nonumber\\
& &- \rho\left( \frac{\xi_x-\eta_xG(\eta)\xi}{1+\eta_x^2}+\varepsilon\varphi_x\right)
\frac{\left(-\eta_x^*G(\eta)\xi-\eta_xG_{\eta}(\eta^*)\xi\right)\left(1+\eta_x^2\right)-2\left(\xi_x-\eta_xG(\eta)\xi\right)\eta_x\eta_x^*}{\left(1+\eta_x^2\right)^2} \nonumber\\
& &+(\bar{\rho}-\rho)g\eta^*+\sigma\frac{\left(1+\eta_x^2\right)^{\frac{3}{2}}\eta_{xx}^*-3\eta_{xx}\eta_x\eta_x^*\left(1+\eta_x^2\right)^{\frac{1}{2}}}{\left(1+\eta_x^2\right)^{3}},
\end{eqnarray}
\begin{eqnarray}\label{eq4.9}
F_{1\bar{\xi}}\bar{\xi}^*& =&c\bar{\rho}\frac{\bar{\xi}^{*}_{x}\eta_x+\bar{G}(\eta)\bar{\xi}^*}{1+\eta^2_x}+ \bar{\rho}\left(\frac{\bar{G}(\eta)\bar{\xi}^*+\eta_x\bar{\xi}^*_x}{1+\eta_x^2}+\varepsilon\bar{\varphi}_y \right)\frac{\bar{G}(\eta)\bar{\xi}^*+\eta_x\bar{\xi}^*_x}{1+\eta_x^2} \nonumber\\
& &+\bar{\rho}\left( \frac{\bar{\xi}^*_x-\eta_x\bar{G}(\eta)\bar{\xi}^*}{1+\eta_x^2}+\varepsilon\bar{\varphi}_x\right)
\frac{\bar{\xi}^*_x-\eta_x\bar{G}(\eta)\bar{\xi}^*}{1+\eta_x^2},
\end{eqnarray}
\begin{eqnarray}\label{eq4.10}
F_{1\xi}\xi^*& =&-c\rho\frac{\xi^{*}_{x}\eta_x+G(\eta)\xi^*}{1+\eta^2_x}- \rho\left(\frac{G(\eta)\xi^*+\eta_x\xi^*_x}{1+\eta_x^2}+\varepsilon\varphi_y \right)\frac{G(\eta)\xi^*+\eta_x\xi^*_x}{1+\eta_x^2} \nonumber\\
& &-\rho\left( \frac{\xi^*_x-\eta_xG(\eta)\xi^*}{1+\eta_x^2}+\varepsilon\varphi_x\right)
\frac{\xi^*_x-\eta_xG(\eta)\xi^*}{1+\eta_x^2},
\end{eqnarray}
\begin{eqnarray}\label{eq4.11}
F_{1c}c^*=\left(\bar{\rho}\frac{\bar{\xi}_{x}\eta_x+\bar{G}(\eta)\bar{\xi}}{1+\eta_x^2}+
\varepsilon\bar{\rho}\bar{\varphi}_{y}-\rho\frac{\xi_{x}\eta_x+G(\eta)\xi}{1+\eta_x^2}-\varepsilon\rho\varphi_{y}\right)c^*;
\end{eqnarray}
\begin{eqnarray}\label{eq4.12}
F_{2\eta}\eta^*=c\eta^*,\quad F_{2\bar{\xi}}\bar{\xi}^*=\bar{\xi}^*,\quad F_{2\xi}\xi^*=0,\quad F_{2c}c^*=\eta c^*;
\end{eqnarray}
\begin{eqnarray}\label{eq4.13}
F_{3\eta}\eta^*=c\eta^*,\quad F_{3\bar{\xi}}\bar{\xi}^*=0,\quad F_{2\xi}\xi^*=\xi^*,\quad F_{2c}c^*=\eta c^*;
\end{eqnarray}
\begin{equation}\label{eq4.14}
F_{4\eta}\eta^*=(H_{\eta}(\eta^*)\xi)_y(\mathbf{z}),\quad F_{4\bar{\xi}}\bar{\xi}^*=0,\quad F_{4\xi}\xi^*=(H(\eta)\xi^*)_y(\mathbf{z}),\quad F_{4c}c^*=c^*.
\end{equation}
Based on the formation of (\ref{eq4.8})-(\ref{eq4.14}), it is enough for us to estimate $\|\eta^*\|_{H^s}$. By (\ref{eq4.8}), it is easy to see that
$$
\eta_{xx}^*=\frac{1}{\sigma}\left(1+\eta_x^2\right)^{\frac{3}{2}}  F_{1\eta}\eta^*-M(\eta^*,\eta^*_x,\eta,\bar{\xi},\xi,c).
$$
We may treat this as an elliptic problem and choose
$$Z^s:=H^{s-1}_{even}(-L,L)\times H^{s-1}_{even}(-L,L)\times H^{s-1}_{even}(-L,L)\times \mathbb{R},$$
which gives
\begin{eqnarray}
\|\eta^*\|_{H^s}& \leq& \|\eta^*\|_{L^2}+\|F_{1\eta}\eta^*\|_{H^{s-2}}+\|M(\eta^*,\eta^*_x,\eta,\bar{\xi},\xi,c)\|_{H^{s-2}}\nonumber\\
&\leq &\|F_{1\eta}\eta^*\|_{H^{s-2}}+\|(\eta^*,\bar{\xi}^*,\xi^*, c^*)\|_{Z^s}.
\end{eqnarray}
Then by using some standard regularity of harmonic and Dirichlet-Neumann operators \cite[Lemma A.1]{SWZ} (see also \cite{Sulem}), we can deduce the estimate (\ref{eq4.7}).

\qed
\section*{Acknowledgments}
G. Dai was supported by National Natural Science Foundation of China (No.
12371110). Y. Zhang was supported by National Natural Science Foundation of
China (No. 12301133), the Postdoctoral Science Foundation of China (No. 2023M741441) and Jiangsu Education Department (No. 23KJB110007).

\section*{Data Availability Statements}
Data sharing not applicable to this article as no datasets were generated or analysed during the current study.

\section*{Conflict of interest}
The authors declare that they have no conflict of interest.

\end{document}